ливість здійснення короткострокового прогнозування поведінки динамічної системи даного технологічного об'єкта.

Збільшення глибини прогнозу можливе за рахунок дослідження і аналізу інших параметрів динамічної системи.

# IDENTIFICATION OF TARGET SYSTEM OPERATIONS. DETERMINATION OF THE TIME OF THE ACTUAL COMPLETION OF THE TARGET OPERATION

**I. Lutsenko**
PhD, Professor
Department of Electronic Devices
Kremenchuk Mykhailo Ostrohradshyi National University
Pervomaiskaya str., 20, Kremenchuk, Ukraine, 39600
E-mail: delo-do@i.ua

*Встановлено, що момент фізичного завершення технологічної операції досліджуваної системи не є моментом завершення цільової операції. Введено поняття моменту фактичного завершення цільової операції. Запропонована система умовних позначень для опису системних процесів цільових операцій. Отримано вирази для чисельного і аналітичного визначення моменту фактичного завершення цільової операції*

*Ключові слова: дослідження операцій, модель цільової операції, момент фактичного завершення цільової операції*

*Установлено, что момент физического завершения технологической операции исследуемой системы не является моментом завершения целевой операции. Введено понятие момента фактического завершения целевой операции. Предложена система условных обозначений для описания системных процессов целевых операций. Получены выражения для численного и аналитического определения момента фактического завершения целевой операции*

*Ключевые слова: исследование операций, модель целевой операции, момент фактического завершения целевой операции*

## 1. Introduction

The process of achieving the goal of any controlled system is carried out by the planning and implementation of individual or linked operations (processes). System operations research is the most delicate instrument, the results of which can be used in solving problems of process optimization, operations planning and analysis. Currently, there is quite a poor set of indicators, allowing to identify the individual system operation. This is caused by the fact that initially the operations research tools were developed for estimating not individual system operations, but the processes of the controlled system within a specified time interval.

Primarily this was due to the accounting system. Subsequently accounting indicators (profit and return) were picked up by economists, for lack of a better alternative, and then by management specialists. System operation, performed within the company (closed system), has no explicit profit, and, consequently, return. In this case we can speak about the added value and the conditional return. Thus, we can say that today there are three basic target operation indicators: the integral value of the input products of operation in comparable cost values – the economic cost ($RE$), the integral value of the output products of operation in comparable cost values – economic income ($PE$), operation time ($T_{op}$) [1]. Derivatives of the basic indicators are the added value (cost) ($PE-RE$) and conditional return ($(PE-RE)/RE$).

## 2. Analysis of literature data and problem statement

Analysis of publications shows that even such a small range of indicators is almost never used by operations research specialists. The main indicator of the operation - efficiency, is not obtained by theoretical investigations based on the key indicators $RE$, $PE$ and $T_{op}$. Most likely, this is why the added value and return are considered economic and not cybernetic indicators. Separate use of these or derived indicators, in general, does not give a clear answer in solving the problems, related to optimal control.

Therefore, specific issues, related to the physics of the processes of the system under investigation [2], problems, where the probability theory [3] or mathematical statistics may be useful, Markov processes [4], in cases when you can directly use the economic indicators [5], etc. are mostly solved in works on the operations research.

There are many works, where performance indicator of operation or process are tried to be replaced with technical indicators [6], probabilistic [7], integration of technological indicators and weighting coefficients [8], by integrating the technological and economic indicators or expert estimates [9].

Works, associated with the target operations research are also hampered by the lack of established notations to describe system processes of target operations.

## 3. Goal and objectives of the paper

The goal of the paper is to obtain a system of target operation indicators, which provide unique identification of the system process with the ability to solve structural and parametric optimization problems in the framework of the controlled system.

Objectives, the solution of which is necessary to achieve the goal are:

– delimitation of the target operation, including in relation to its physical boundaries; numerical and analytical determination of the time of actual completion of the target operation;

– determination of the complex cost of operation on its deployed model;

– numerical and analytical determination of the complex cost of the target operation;

– determination of the potential effect of the target operation;

– determination of the efficiency of use of the target operation resources using numerical methods and analytically.

## 4. The system of notations in the identification problems

Unlike physical models of technological operations, which describe the features of the technological process of converting raw and energy products, models of target operations have no conceptual differences. This means that any target operation can be described based on the unified system of notations. Notations, used in the paper are given below:

$r_i$ – the i-th input product of the system under investigation;

$p_j$ – the j-th output product of the system under investigation;

$rq_i(t)$ – registration signal, displaying quantitative parameter of the i-th input product of the system under investigation;

$pq_j(t)$ – registration signal, displaying quantitative parameter of the j-th output product of the system under investigation;

$rs_i$ – cost estimate of the i-th input product of the target operation;

$ps_j$ – cost estimate of the j-th input product of the target operation;

$re_i(t)$ – registration signal of the i-th input product of the system under investigation, reduced to comparable cost values $re_i(t) = rs_i \int_0^t rq_i(t)dt$;

$pe_j(t)$ – registration signal of the j-th output product of the system under investigation, reduced to comparable cost values $pe_j(t) = ps_j \int_0^t pq_j(t)dt$;

$re(t)$ – total registration signal of the i-input, reduced to comparable cost values, products of the system under investigation $re(t) = \sum_{i=1}^{I} rs_i \int_0^t rq_i(t)dt$;

$pe(t)$ – total registration signal of the j-output, reduced to comparable cost values, products of the system under investigation $pe(t) = \sum_{j=1}^{J} ps_j \int_0^t pq_j(t)dt$;

$t_s$ – time of the beginning of the target operation;
$t_f$ – time of the physical completion of the target operation;
$t_a$ – time of the actual completion of the target operation;
$t_r$ – time of the supply of input products of the target, reduced, simplified operation;
$t_p$ – time of the transfer of output products of the target, reduced, simplified operation;

$RE = \int_{t_0}^{t_a} re(t)dt$ – cost estimate of input products of the target operation;

$PE = \int_{t_0}^{t_a} pe(t)dt$ – cost estimate of output products of the target operation;

$ire(t) = \int_0^t re(t)dt$ – thread of resource consumption of the target operation;

$ipe(t) \int pe(t)dt$ – thread of resource productivity of the target operation;

$ice(t) = \int_0^t re(t)dt + \int_0^t pe(t)dt$ – single-threaded model of the target operation;

$vre(t) = \int_0^t ire(t)dt$ – integral function of the thread of resource consumption;

$vpe(t) = \int_0^t ipe(t)dt$ – integral function of thread of resource productivity.

## 5. Determination of the time of the physical completion of the target operation

Determination of the time of completion of a separate system operation is considered a fairly simple matter. System operation is considered complete when the whole raw or buffered product, fed to the system input goes beyond the system. It is most convenient to analyze the time frames of the operation model on the example of the movement of one product through the buffering system. Energy costs are neglected.

Let us assume that input portion of the product $r$ is fed to the buffering system input. After some time, the output portion of the product $p$ has gone beyond the buffering system (Fig. 1).

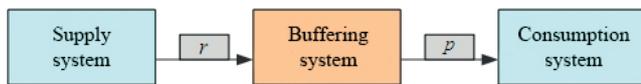

Fig. 1. Movement of the product through the buffering system

Since the buffering system does not alter the product properties, the volume of the input product $r$ is equal to the volume of the output product $p$. The time $t_s$ is the time of the beginning of the supply of the input product and the time $t_f$ – the time of the completion of the feeding of the output product (Fig. 2).

If registration of the movement of the input product $r$ at the system input can be described by function $rq(t)$, and registration of the movement of the product $p$ at the system output - function $pq(t)$, the function of internal reserves may be determined by the expression $icq(t) = \int_0^t rq(t)dt + \int_0^t pq(t)dt$.

Time $t_s$ is considered the time of the beginning of the system operation, and time $t_f$ – the time of its completion.

Time $t_f$ is characterized by the fact that the whole internal product of the system under investigation has been fully transferred on its output by this time. This time is defined in the paper as **the time of the physical completion of the system technological operation**.

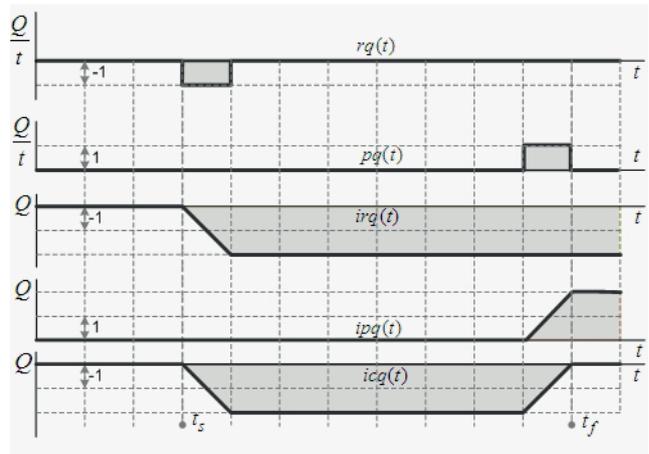

Fig. 2. Technology of determining the level of reserves in the buffering system

However, the goal of any system operation is not a conversion, transformation or displacement of the input system products for the preparation of the final product, but increase in the value of the output system products with respect to the value of its input products.

If we define the cost estimate of the unit of the input product as the symbols $rs$, and units of the output product as $ps$, the deployed model of the products movement in comparable values will be of the form

$ice(t) = rs\int_0^t rq(t)dt + ps\int_0^t pq(t)dt = ire(t) + ipe(t)$ [1] (Fig. 3).

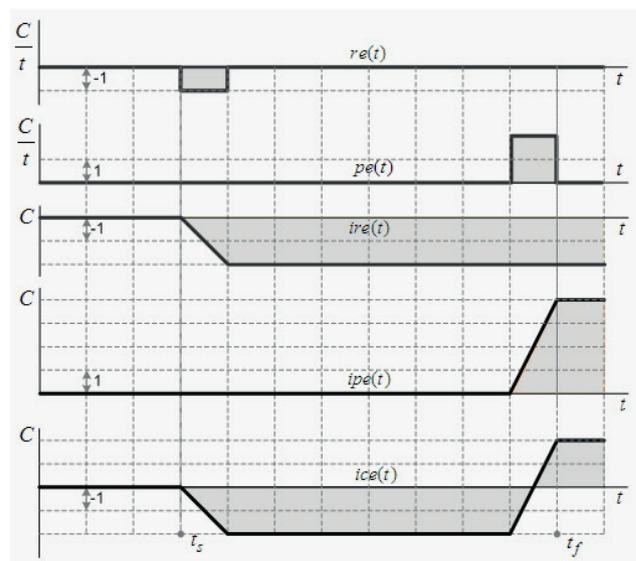

Fig. 3. Target operation of the buffering system

As seen, the appearance of the function $ice(t)$ is quite different from the function $icq(t)$. Are these differences formal or it is something more?

## 6. Determination of the time of actual completion of the target operation

Let us proceed from the deployed model of the operation (Fig. 3) to a simplified reduced deployed model of the operation [1] (Fig. 4) and pay attention to the closed thread of tight resources $ibe(t)$. In the time interval from the time $t_r$ until the time $t_p$, input technological products of the operation were linked by internal processes of the system under investigation. The function $vbe(t) = \int_0^t |ibe(t)| dt$ quantitatively displays in time the estimated – time component, related to tight resources. On the other hand, the function $vde(t) = \int_0^t ide(t) dt$ quantitatively displays in time the estimated – time component, related to the released resources of the operation.

The intersection point of these two functions indicates the time when the target thread $ide(t)$ compensates in time the thread of tight resources. Let us define this time as the time of the actual completion of the operation (TACO). This time is defined on the graph as $t_a$.

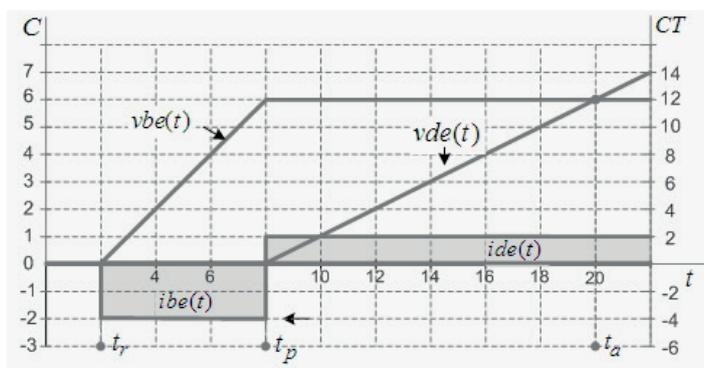

Fig. 4. Determination of the time of the logical completion of the operation with the use of integral functions of the thread of tight resources and the target thread

As can be seen, the process of conversion of input products requires a certain time. For the presented operation model, this time is 6 time intervals. At the time of physical completion of the operation ($t_f = 8$), the process of conversion of input products is completed.

Analyzing the deployed model of the operation, we can see that its implementation is accompanied by losses, caused by binding of input products in time. These losses, displayed by the closed thread of tight resources $ibe(t)$ are the greater the higher the cost estimate of tight products and the longer the process of their conversion. Consequently, the geometry of the closed thread of the negative half-plane displays the system losses, associated with the achievement of the operation objective.

At the time of the physical completion of the operation ($t_f$), system under investigation transmits the output product to the consumption system. This means that the value of the input products is increased, and the output product is transferred to gain the effect, resulting from its use.

In efficient operation, at the time of its physical completion, the maximum value of the thread $ipe(t)$ is above the maximum value of the thread $ire(t)$. Due to this difference in the maximum values, the target thread $ide(t)$ is formed.

The value of the thread $ide(t)$ is caused by the difference in the cost estimate of output products of the operation and input products of the operation. But the thread of tight resources $ibe(t)$ can not be compensated by the cost estimate of the output product. This is evident during formal comparison of measurement units. Thread, as the area has the measurement unit «*Cost* x *Time*», and the measurement of the output product *p* is determined by the category of «*Cost*». Measurement units of these categories are different, therefore, they are incompatible.

You can compare, for example, the thread of tight resources $ibe(t)$ and the target thread $ide(t)$. It is logical to assume that the operation can not be considered complete until the target thread $ide(t)$ compensates for the thread of resources $ibe(t)$, tied by the operation.

For the operation under investigation, this time is $t_a$. At this time, the area of the thread $ide(t)$ becomes equal to the area of the thread $ibe(t)$.

We succeeded to determine the TACO in such a simple way only because registration signals were selected for the studied simplified reduced operation so that the areas of threads can be easily compared. To determine the TACO of any effective operation, threads $ibe(t)$ and $ide(t)$ must be integrated. Let us denote integral function of the module of the thread of tight resources as $vbe(t)$ and integral function of the target thread $ide(t)$ as $vde(t)$ (Fig. 4). Then

$$vbe(t) = \int_{t_0}^{t} |ibe(t)| dt; \quad vde(t) = \int_{t_0}^{t} ide(t) dt. \quad (1)$$

The measurement unit of integral functions $vde(t)$ and $vbe(t)$ is defined by the product of the expert (cost) component by the time component. This measurement unit will be denoted as CT.

Thus, the graphics of functions $vde(t)$ and $vbe(t)$ intersect at point of equality of areas of threads $ide(t)$ and $ibe(t)$.

The time of the actual completion of the operation can be defined as *the time of equality of integral functions from the function of the thread of tight resources of the target thread.*

## 7. Analytical determination of the time of the actual completion of the operation

Determination of the time of the physical completion of the operation (TPCO) based on the thread of tight resources $ibe(t)$ and target thread $vde(t)$ allows to understand one of the process aspects in achieving the goal of the operation. However, TPCO can be defined in an easier way based on the threads of resource consumption $ire(t)$ and resource productivity $ipe(t)$ directly, without the intermediate conversion of these threads into the threads $ibe(t)$ and $ide(t)$.

Let us construct the integral function of the module of the thread $ire(t)$, which we denote as $vre(t)$ and integral function of the thread $ipe(t)$, which we denote as $vpe(t)$ (Fig. 5).

Functions $vre(t)$ and $vpe(t)$ intersect at the point, which we defined as the TPCO. As can be seen, determination of the TPCO with the use of the resource consumption and resource productivity threads requires less computational operations.

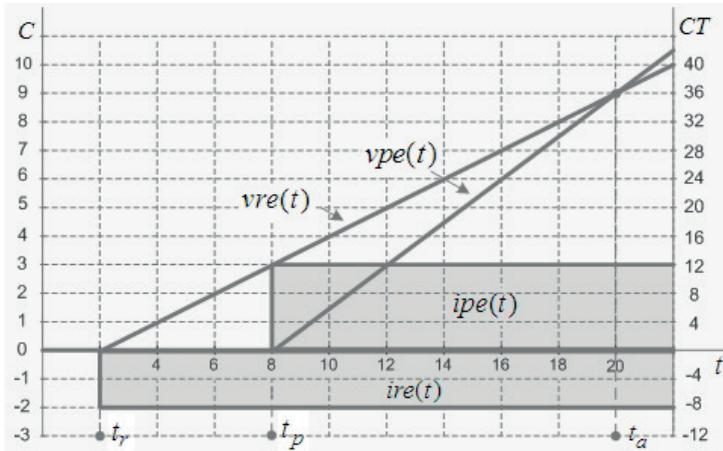

Fig. 5. Determination of the time of the actual completion of the operation with the use of integral functions of the resource consumption and resource productivity threads

Thus, TPCO can also be defined as the time of *equality of integral functions from the functions of the resource consumption and resource productivity threads.* But the fact that the nonlinear functions vre(t) and vpe(t), unlike the essentially non-linear functions vbe(t) and vde(t) can be replaced by linear functions $vre^*(t)$ and $vpe^*(t)$ is much more important (Fig. 6).

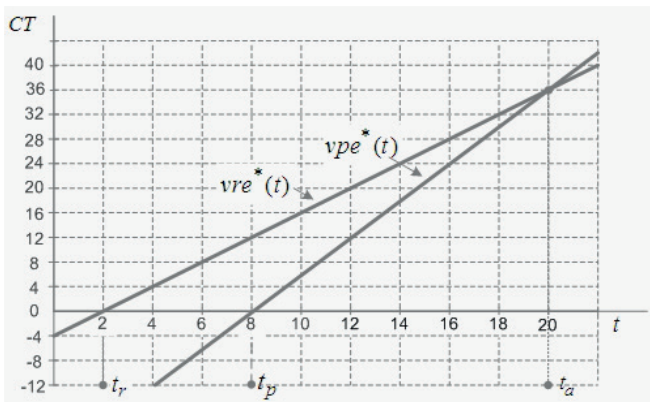

Fig. 6. Determination of the time of the actual completion of the operation with the use of linear functions $vre^*(t)$ and $vpe^*(t)$

$$vre^*(t) = ire^*(t) \cdot t - C_r + C, \tag{2}$$

$$vpe^*(t) = ipe^*(t) \cdot t - C_p + C. \tag{3}$$

The constants $C_r$ and $C_p$ respectively, are defined by the expressions

$$C_r = ire^*[t_a] \cdot t_a, \tag{4}$$

$$C_p = ipe^*[t_a] \cdot t_a. \tag{5}$$

Notation, for example $ipe^*[t_a]$, means that the value of the function $ipe^*(t)$ at a point $t_a$ is taken. Constant C determines the displacement of functions $vre^*(t)$ and $vpe^*(t)$. For functions $vre^*(t)$ and $vpe^*(t)$ (Fig. 6) C = 36.

With regard to (4) and (5), expressions (2) and (3) can be written as

$$vre^*(t) = ire^*(t) \cdot t - ire^*[t_a] \cdot t_a + C, \tag{6}$$

$$vpe^*(t) = ipe^*(t) \cdot t - ipe^*[t_a] \cdot t_a + C. \tag{7}$$

From the expressions

$$0 = ire^*(t) \cdot t - ire^*[t_a] \cdot t_a + C, \tag{8}$$

$$0 = ipe^*(t) \cdot t - ipe^*[t_a] \cdot t_a + C. \tag{9}$$

we obtain the values of the time points, when the functions $vre^*(t)$ and $vpe^*(t)$ cross it

$$t_r = \left(ire^*[t_a] \cdot t_a - C\right) / ire^*[t_r], \tag{10}$$

$$t_p = \left(ipe^*[t_a] \cdot t_a - C\right) / ipe^*[t_p]. \tag{11}$$

where $t_r$ – time of crossing of the time axis by the function $vre^*(t)$; $t_p$ – time of crossing of the time axis by the function $vpe^*(t)$.

It is possible to check that the time $t_r$ corresponds to the time of the registration of input products of the focused operation, time $t_p$ corresponds to the time of the registration of output products of the focused operation.

Let us write the system of equations

$$\begin{cases} 0 = ire^*[t_r] \cdot t_r - ire^*[t_a] \cdot t_a, \\ 0 = ipe^*[t_p] \cdot t_p - ipe^*[t_a] \cdot t_a. \end{cases} \tag{12}$$

Having solved it relatively to $t_a$, we obtain

$$t_a = \frac{ipe^*[t_p] \cdot t_p - ire^*[t_r] \cdot t_r}{ipe^*[t_p] - ire^*[t_r]}. \tag{13}$$

Given that for simplified reduced operations, $ipe^*[t_p]$ is numerically equal to PE, and $ire^*[t_r]$ is numerically equal to the value $|RE|$, the numerical value of $t_a$ can be determined using the following expression

$$t_a = \frac{PE \cdot t_p - |RE| \cdot t_r}{PE - |RE|}. \tag{14}$$

For example, for the operation under investigation, we obtain

$$t_a = \frac{PE \cdot t_p - |RE| \cdot t_r}{PE - |RE|} = \frac{3 \cdot 8 - 2 \cdot 2}{3 - 2} = 20.$$

Equation (14) can be written in general form to determine the TPCO of any effective operation

$$t_a = \frac{\int_0^{t_f}\left[pe(t)\cdot t\right]dt - \int_0^{t_f}\left[|re(t)|\cdot t\right]dt}{\int_0^{t_f} pe(t)dt - \int_0^{t_f}|re(t)|dt}. \qquad (15)$$

In discrete coordinate systems, expression for determining the TPCO takes the form

$$n_a = \frac{\sum_{k=1}^{K}\left[pe_j n_j\right] - \sum_{k=1}^{K}\left[re_i n_i\right]}{\sum_{k=1}^{K} pe_j - \sum_{k=1}^{K} re_i},\ i=\overline{1,K};\ j=\overline{1,K}. \qquad (16)$$

Examples of practical application of the formulas (14) and (16) are available in the resource [10].

## 8. Conclusions

It was found that the time of the physical completion of the technological operation of the system under investigation is not the time of the completion of the target operation. The concept of the time of the actual completion of target operation was introduced. A system of notations to describe the system processes of target operations was proposed.

It was found that the time of completion of system operation with regard to achieving the goal is defined by equality of the integral estimates of thread of tight resources of the operation and the target thread. The expressions for the numerical and analytical determination of the time of actual completion of the target operation were obtained.